\documentclass{amsart}

\usepackage{amscd,color,tikz,latexsym,amssymb}
\usetikzlibrary{calc}
\usetikzlibrary{arrows}
\usepackage{tikz-cd}

\usepackage{amssymb,tikz}
\usepackage{amscd}
\usepackage{latexsym}
\usepackage{amsfonts}
\usepackage{amsmath}
\usepackage[latin1]{inputenc}

\newtheorem{proposition}{Proposition}
\newtheorem{theorem}{Theorem}
\newtheorem{definition}{Definition}
\newtheorem{remark}{Remark}
\newtheorem{conjecture}{Conjecture}

\newcommand{\cla}{c_\la}

\newcommand{\IMUN}{\mathsf i}
\newcommand{\imun}{\mathsf i}

\newcommand{\INTEGERS}{\mathbb Z}

\newcommand{\la}{\mathsf{b}}

\newcommand{\QDILOG}{\Phi_{\la}}

\newcommand{\REALS}{\mathbb R}

\newcommand{\im}{\mathop{\fam0 Im}\nolimits}

\newcommand{\Hom}{\mathop{\fam0 Hom}\nolimits}

\newcommand{\vol}{\mathop{\fam0 Vol}\nolimits}

\newcommand{\sign}{\mathop{\fam0 sign}\nolimits}

\newcommand{\bC}{{\mathbb C}}
\newcommand{\bR}{{\mathbb R}}

\newcommand{\Z}{{\mathbb Z}}
\newcommand{\bZ}{\Z{}}

\newcommand{\D}{{\mathcal D}}
\newcommand{\ra}{\mathop{\fam0 \rightarrow}\nolimits}

\renewcommand{\L}{{\mathcal L}}

\newcommand{\SU}{\mathop{\fam0 SU}\nolimits}

\newcommand{\CC}{\mathcal{B}}

\renewcommand{\S}{\mathcal{S}}

\newcommand{\wf}{\mathrm{WF}}

\newcommand{\parfun}{\mathop{Z_\hbar}}

\newcommand{\BQDILOG}{\bar{\Phi}_{\la}}
\newcommand{\Imun}{\mathsf{i}}

\newcommand{\bS}{\mathbb{S}}

\begin{document}

\title{The Teichm\"{u}ller TQFT}

\author{J{\o}rgen Ellegaard Andersen}
\address{Center for Quantum Geometry of Moduli Spaces\\
Department of Mathematics\\
        University of Aarhus\\
        DK-8000, Denmark}
\email{andersen@qgm.au.dk}

\author{Rinat Kashaev}
\address{University of Geneva\\
2-4 rue du Li\`evre, Case postale 64\\
 1211 Gen\`eve 4, Switzerland}
\email{rinat.kashaev@unige.ch}

\thanks{Supported in part by the center of excellence grant ``Center for quantum geometry of Moduli Spaces" DNRF95, from the Danish National Research Foundation, Swiss National Science Foundation and the ESF-ITGP network.}

\begin{abstract}
We review our construction of the Teichm\"{u}ller TQFT. We recall our volume conjecture for this TQFT and the examples for which this conjecture has been established. We end the paper with a brief review of our new formulation of the Teichm\"{u}ller TQFT together with some anticipated future developments.
\end{abstract}

\maketitle

\section{Introduction}\label{intro}

Topological Quantum Field Theories (TQFT's) were discovered and axiomatised by Atiyah~\cite{At}, Segal~\cite{S} and Witten~\cite{W}. Following Witten's suggestions in \cite{W}, the first examples in $2+1$ dimensions were constructed by Reshetikhin and Turaev \cite{RT1,RT2,T} based on the representation theory of quantum groups at roots of unity. The resulting {\em Witten--Reshetikhin--Turaev} TQFT (WRT-TQFT) has also been formulated in pure topological terms in \cite{BHMV1,BHMV2} and it was conjectured by Witten in \cite{W} to be related to quantum conformal field theory and geometric quantization of moduli spaces. It has been further developed in \cite{TUY,ADW,H} and in a series of papers including the work of Laszlo~ \cite{La1} who proved that the Hitchin and the TUY connections agree in the closed surface case. The equivalence of the geometric and combinatorial constructions has been finally verified in \cite{AU1,AU2,AU3,AU4} by the first author of this paper jointly with Ueno and exploited in \cite{A5, A,A2, A3, AHi, AJ1, AHJMM} to establish some strong properties of the WRT-TQFT. 

In parallel to the surgery based construction of the WRT-TQFT by Reshetikhin and Tureav, there is the Turaev--Viro construction of the TV-TQFT~ \cite{TV}, which is also a combinatorial construction, but it uses triangulations instead, and where Reshetikhin and Turaev had to prove the invariance under the surgery presentation of their construction, Turaev and Viro proved invariance under the Pachner 2-3, 3-2, 4-1 and 1-4 moves of theirs. It turned out that the TV-TQFT was the Hermitian endomorphism theory of the WRT-TQFT.

In his paper \cite{W2}, Witten further proposed that quantum Chern--Simon theory for non-compact groups should also exist as generalised TQFT's with underlying infinite dimensional state vector spaces. A series of papers on this subject have subsequently emerged in the physics literature, including \cite{MR1133274, D,DGG,DGLZ,DG,GM,G,Hik1,Hik2,W3}. However, a mathematical definition of these theories has been lacking for a long time.

In a series of papers \cite{AK1,AK2,AK3,AK4,AK5}, the authors of this paper, have provided a rigorous construction of such a TQFT, known as the \emph{Teichm\"{u}ller TQFT}.  Our construction uses combinatorics of $\Delta$-complexes with fixed number of vertices which we call triangulations and it builds on quantum Teichm\"{u}ller theory, as developed by Kashaev \cite{K1}, and Chekhov and Fock \cite{CF}, which produces unitary representations of centrally extended mappings class groups of punctured surfaces in infinite-dimensional Hilbert spaces. In this paper we shall first review our original formulation presented in \cite{AK1,AK2,AK3}. 


The central ingredients in quantum Teichm\"{u}ller theory are, on the one hand, Penner's coordinates of the decorated Teichm\"{u}ller space and the Ptolemy groupoid \cite{Pen1} with applications summarised in \cite{Pen2} and, on the other hand, Faddeev's quantum dilogarithm \cite{F} which finds its origins and applications in quantum integrable systems \cite{FKV,BMS1,BMS2,Te}. Faddeev's quantum dilogarithm has already been used  in formal state-integral constructions of perturbative invariants of three manifolds in the works \cite{Hik1,Hik2,DGLZ,DFM,D}, but without addressing the  important questions of convergence or triangulation independence.

There are further ingredients which we had to introduce in \cite{AK1} in order to lift quantum Teichm\"{u}ller theory to a TQFT. The important one is the weight function for tetrahedra, whose edges are labeled by dihedral angles of hyperbolic ideal tetrahedra. 
It is not immediately clear what are the topological invariance properties of our TQFT which depends on those dihedral angles. It turns out,  however, the partition function of a given triangulation is invariant under certain Hamiltonian gauge group action in the space of angles so that the corresponding symplectically reduced space is determined by the total dihedral angles around edges and the first cohomology group of the (cusp) boundary. As a consequence, under the condition that the triangulation in question is such that the second homology group of the complement of the vertices is trivial, the partition function descends to a well defined function on an open convex subset of this reduced angle space (corresponding to strictly positive angles). Furthermore, if we have two triangulations admitting angle structures (which correspond to balanced edges with total dihedral angles equal to $2\pi$) and related by a Pachner 2-3 or 3-2 move, then the two convex subsets intersect non-trivially and the two partitions functions agree on the overlap. The additional fact that the partition functions depend analytically on the dihedral angles implies that their common restriction to the overlap completely determines both of them and it is in this sense that our TQFT is topologically invariant. 


The partition functions of our TQFT take their values in the vector spaces of tempered distributions over euclidian spaces which do not form a category, since it is not always possible to multiply and push forward tempered distributions. Instead, they form what we call a categroid that is the same as a category, except that we are allowed to compose not all morphisms which are composable in the categorical sense, but only a subset thereof (which we review in Section \ref{target}). Symmetrically, the domain of our TQFT, the set of oriented triangulated pseudo $3$-manifolds, also forms only a categroid, due to the above mentioned homological condition on triangulations.


We shall further review a version of the volume conjecture for the Teichm\"{u}ller TQFT, which states that the partition function {\em decays} exponentially fast in Planck's constant  with the rate given by the hyperbolic volume of the manifold. 

Interestingly, due to subsequent developments, we have now at least two formulations of the Teichm\"{u}ller TQFT. The original formulation, which is defined only for admissible pseudo 3-manifolds (see Definition~\ref{admis} below), and the new formulation, which does not impose any restrictions on the topology of pseudo $3$-manifolds. We will briefly discuss the new formulation in the end of this paper, together with a number of future developments which we anticipate.

The paper is organised as follows. In Section~\ref{domain}, we review the domain categroid, on which our original formulation of the Teichm\"{u}ller TQFT is defined, while the target categroid is reviewed in Section~\ref{target}. In Section~\ref{functor}, we review our TQFT functor between these two categroids and state the main Theorem~\ref{MT}, proved in \cite{AK1}, which establishes the well definedness of the functor. In Section~\ref{vol}, we formulate the volume conjecture for the Teichm\"{u}ller TQFT and  describe a couple of examples for which that conjecture has already been established. In the final Section~\ref{future}, we briefly describe the new formulation of the Teichm\"{u}ller TQFT and anticipate a number of future developments.

\subsection*{Acknowledgements} We would like to thank Gregor Masbaum, Nikolai Reshetikhin, and Vladimir Turaev  for valuable discussions. Our special thanks go to Feng Luo for  explaining to us the topological significance of non-negative angle structures and Anton Mellit for rising important questions about analytically continued partition functions.

\section{The topological domain categroid}
\label{domain}


In this section we set up the topological domain categroid on which our Teichm\"{u}ller TQFT is defined. Since, in its original formulation, this TQFT is distributional in nature, we cannot use the full 3-dimensional bordism category of triangulated pseudo 3-manifolds (with extra structures), and we consider a suitable sub-categroid of it as follows.

Let $X$ be a triangulated pseudo $3$-manifold, that is a CW-complex obtained by gluing finitely many tetrahedra with ordered vertices along codimension one faces with respect to order preserving simplicial maps.

For $i\in\{0,1,2,3\}$, we will denote by $\Delta_i(X)$ the set of $i$-dimensional cells in $X$. For any $i>j$, we also denote

\[
\Delta_{i}^{j}(X):=\{(a,b)\mid a\in\Delta_i(X),\ b\in\Delta_j(a)\},
\]
where the cell $a$, when considered as a CW-complex, is taken in the form of the standard simplex without identifications on its boundary induced by gluings.

\begin{figure}[h]
\centering
\includegraphics{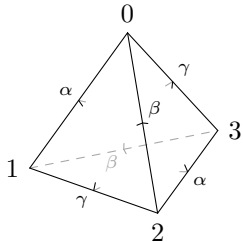}
\caption{A tetrahedron with ordered vertices and dihedral angles.}
\label{fig:Diangles}
\end{figure}

A \emph{shape structure} on $X$  is an assignment  to each edge of each tetrahedron of $X$ a positive number,
$
\alpha_X\colon \Delta_{3}^{1}(X)\to\bR_{>0},
$
 called \emph{dihedral angles} such that the sum of the angles at any three edges sharing a vertex of a tetrahedron is $\pi$.  It is straightforward to see that the dihedral angles at opposite edges of all tetrahedra are equal, so that each tetrahedron acquires three dihedral angles associated to the three pairs of opposite edges which sum up to $\pi$, see Fig.~\ref{fig:Diangles}.  In other words, a shape structure provides each tetrahedron with the geometric structure of an ideal hyperbolic tetrahedron. An oriented triangulated pseudo 3-manifold with a shape structure is called a \emph{shaped pseudo 3-manifold}. We denote the set of shape structures on $X$ by $S(X)$.

An edge is called \emph{balanced} if it is internal and the sum of dihedral angles around it is $2\pi$. An edge which is not balanced is called \emph{unbalanced}. An \emph{angle structure} on a closed triangulated pseudo $3$-manifold, introduced by Casson, Rivin and Lackenby \cite{C,R,Lackenby}, is a shape structure where all edges are balanced.

We will also consider the situation, where we are given a one dimensional sub-complex $\Gamma$ of $\Delta_1(X)$, such that all univalent vertices of $\Gamma$ are on the boundary of $X$ (such a sub-complex we will call an allowed one dimensional sub-complex $\Gamma \subset X$). 

We extend the shape structure by a real parameter called the {\em level}. This is an analog of the framing in the context of the WRT-TQFT.
Thus, a \emph{levelled shaped pseudo 3-manifold} is a pair $(X,\ell_X)$ consisting of a shaped pseudo 3-manifold $X$ and a level  $\ell_X\in\bR$, and we denote by $LS(X)$ the set of all levelled shaped structures on $X$. 
 
 Two levelled shaped pseudo 3-manifolds $(X,\ell_X)$ and $(Y,\ell_Y)$ are called \emph{gauge equivalent} if there exists an isomorphism $h\colon X\to Y$ of the underlying cellular structures and a function
\[
g\colon \Delta_1(X)\to \REALS, \quad g\vert_{\Delta_1(\partial X)}=0,
\]
such that the shape structures of $X$ and $Y$ and levels $\ell_X, \ell_Y$ are related by the formulae stated in Definition 2 of \cite{AK1}.
As is explained in Section~2 of \cite{AK1}, this equivalence is  induced by a Hamiltonian group action corresponding to the Neumann--Zagier symplectic structure.
In the particular case $X=Y$ and the identity isomorphism, we get the notion of \emph{based gauge equivalence} of levelled shaped pseudo $3$-manifolds. The set of based gauge equivalence classes of levelled shape structures on $X$ is denoted $LS_r(X)$ and we denote by $S_r(X)$ the corresponding set of based gauge equivalence classes of just shape structures (obtained by forgetting the level).

By removing the positivity condition in the definition of a shape structure, we define a \emph{generalised shape} structure on $X$, and we denote by $\tilde{S}(X)$ the set of generalised shape structures. Levelled generalised shaped structures as well as their gauge equivalence are defined analogously. The space of based gauge equivalence classes of generalised shape structures (respectively levelled generalised shaped structures) is denoted $\tilde{S}_r(X)$ (respectively $\widetilde{LS}_r(X)$).
Remark that $S_r(X)$ is an open convex subset of $\tilde{S}_r(X)$.

Let 
$
\tilde{\Omega}_X : \tilde{S}(X) \ra \bR^{\Delta_1(X)}
$
be the map which associates to an edge $e$ the sum of the dihedral angles around $e$. The values of $\tilde{\Omega}_X$ will be called \emph{(edge) weights}.  Due to gauge invariance, $\tilde{\Omega}_X$ induces a unique map
$
\tilde{\Omega}_{X,r} : \tilde{S}_r(X) \ra \bR^{\Delta_1(X)}.
$

Let $N_0(X)$ be a sufficiently small tubular neighbourhood of $\Delta_0(X)$. The boundary $\partial N_0(X)$ is a two dimensional surface, which is possibly disconnected and possibly with boundary, if $\partial X \neq \varnothing$. Theorem~1 of \cite{AK1} states that the map
$
\tilde{\Omega}_{X,r}
$
 is an affine $H^1(\partial N_0(X), \bR)$-bundle.
The space $\tilde{S}_r(X)$ carries a Poisson structure whose symplectic leaves are the fibers of $\tilde{\Omega}_{X,r}$ and which is identical to the Poisson structure induced by the $H^1(\partial N_0(X), \bR)$-bundle structure. The natural projection map from $\widetilde{LS}_r(X) $ to $\tilde{S}_r(X)$ is an affine  $\bR$-bundle which restricts to the affine $\bR$-bundle $LS_r(X)$ over $S_r(X)$.

If $h:X \ra Y$ is an isomorphism of cellular structures, then we get an induced Poisson isomorphism $h^* : \tilde{S}_r(Y) \ra \tilde{S}_r(X)$ which is an affine bundle isomorphism with respect to the induced group homomorphism
\[
h^* : H^1(\partial N_0(Y), \bR)\ra H^1(\partial N_0(X), \bR)
 \]
 and which maps $S_r(Y)$ to $S_r(X)$. Furthermore, $h$ induces an isomorphism
\[
h^*\colon \widetilde{LS}_r(Y) \to \widetilde{LS}_r(X)
 \]
 of affine $\bR$-bundles covering the map
 \(
 h^*\colon \tilde{S}_r(Y) \to \tilde{S}_r(X)
 \)
  and which also maps $LS_r(Y)$ to $LS_r(X)$.

\begin{figure}[htbp]
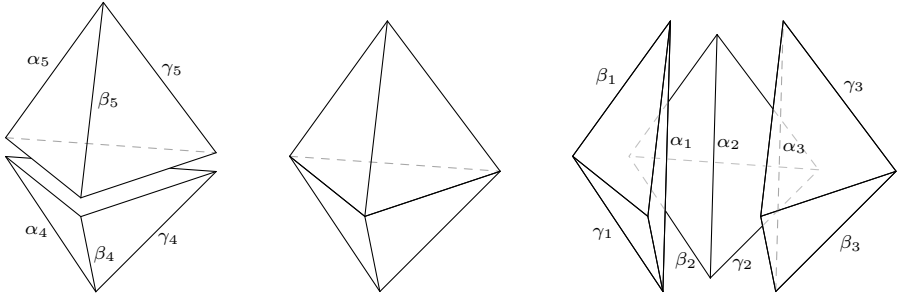

\centering
\includegraphics{tetra-4.mps}
\hskip.8cm
\includegraphics{tetra-5.mps}
\hskip.8cm
\includegraphics{tetra-6.mps}
\caption{The 3-2 Pachner move}
\label{fig:3.2}
\end{figure}

Let us now consider the 3-2 Pachner move illustrated in Fig.~\ref{fig:3.2}.
Let $e$ be a balanced edge of a shaped pseudo 3-manifold $X$ and assume that  $e$ is shared by exactly three distinct tetrahedra $t_1,t_2,t_3$. Let $S$ be a shaped pseudo 3-submanifold of $X$ composed of the tetrahedra $t_1,t_2,t_3$. Note that $S$ has $e$ as its only internal and balanced edge. There exists another triangulation $S_e$ of the topological space underlying $S$ such that the triangulation of $\partial S$ coincides with that of $\partial S_e$, but which consists of only two tetrahedra $t_4,t_5$. We see that this change has the effect of removing the edge $e$ so that $\Delta_1(S_e)=\Delta_1(S)\setminus\{e\}$. Moreover, there exists a unique shaped structure on $S_e$ which induces the same edge weights as the shape structure of $S$. For shape variables $(\alpha_i,\beta_i,\gamma_i)$ for $t_i$ (where $\alpha_i$ are the angles at $e$), the explicit map is given by

\begin{equation}\label{P32a}\begin{array}{ll}
\alpha_4 = \beta_2 + \gamma_1 & \alpha_5 = \beta_1 + \gamma_2\\
\beta_4 = \beta_1 + \gamma_3 & \beta_5 = \beta_3 + \gamma_1\\
\gamma_4 = \beta_3 + \gamma_2 & \gamma_5 = \beta_2 + \gamma_3.
\end{array}
\end{equation}
We observe that the equation $\alpha_1+\alpha_2 + \alpha_3 = 2\pi$ implies that the angles for $t_4$ and $t_5$ sum up to $\pi$. Moreover the positivity of the angles for $t_1, t_2, t_3$ implies that the angles for $t_4$ and $t_5$ are also positive. On the other hand, it is not automatic that we can solve for positive angles for $t_1,t_2,t_3$ given the positive angles for $t_4$ and $t_5$. However if we have two positive solutions for the angles for $t_1,t_2,t_3$ for the same $t_4,t_5$, then they are gauge equivalent and satisfy the equality $\alpha_1+\alpha_2 + \alpha_3 = 2\pi$.

\begin{definition}
We say that a shaped pseudo 3-manifold $Y$ is obtained from $X$ by a \emph{shaped 3-2 Pachner move} along the edge $e$ if $Y$ is obtained from $X$ by replacing $S$ by $S_e$, and we write $Y=X_e$.
\end{definition}

We observe from the above that there is a canonical map $ P^e : S(X) \ra  S(Y),$ which naturally extends to a map $ \tilde{P}^e : \tilde S(X) \ra  \tilde S(Y).$
We get the following commutative diagram
\begin{equation}
\begin{CD}
\tilde \Omega_X(e)^{-1}(2\pi) @>{\tilde P^e}>> \tilde S(Y)\\
@VVV @VVV\\
\tilde \Omega_{X, r} (e)^{-1}(2\pi)  @>\tilde P_r^e>> \tilde S_r(Y)\\
@VV{\operatorname{proj}\ \circ \ \tilde \Omega_{X,r}}V @VV{\tilde \Omega_{Y,r}}V\\
\bR^{\Delta_1(X)\setminus\{e\}} @> = >>  \bR^{\Delta_1(Y)}
\end{CD}
\end{equation}

Moreover
$$\tilde P_r^e( \tilde \Omega_{X, r} (e)^{-1}(2\pi)  \cap S_r(X)) \subset S_r(Y).$$
In particular, we observe that if $\tilde \Omega_{X, r} (e)^{-1}(2\pi)  \cap S_r(X) \neq \emptyset$ then $S_r(Y)\neq \emptyset$, but the converse is not necessarily true. The following theorem is proved in Section~2 of \cite{AK1}.

\begin{theorem}\label{3-2SS}
Suppose that  a shaped pseudo 3-manifold $Y$ is obtained from a shaped pseudo 3-manifold $X$ by a levelled shaped 3-2 Pachner move along the edge $e$.
Then the map $\tilde P_r^e$ is a Poisson isomorphism, which is covered by an affine $\bR$-bundle isomorphism from $\widetilde{LS}_r(X)|_{\tilde \Omega_{X, r} (e)^{-1}(2\pi) }$ to $\widetilde{LS}_r(Y)$.
\end{theorem}

We also say that a levelled shaped pseudo 3-manifold $(Y,\ell_Y)$ is obtained from a levelled shaped pseudo 3-manifold $(X,\ell_X)$ by a \emph{levelled shaped 3-2 Pachner move} if there exists $e\in\Delta_1(X)$ such that $Y=X_e$ and the levels are related by the formula stated just above Definition 9 in \cite{AK1}.

\begin{definition}\label{PMove}
A (levelled) shaped pseudo 3-manifold $X$ is called a \emph{Pachner refinement} of a (levelled) shaped pseudo 3-manifold $Y$ if there exists a finite sequence of (levelled) shaped pseudo 3-manifolds
$
X=X_1,\ X_2,\ \ldots,\ X_n=Y
$
such that for any $i\in\{1,\ldots,n-1\}$, $X_{i+1}$ is obtained from $X_i$ by a (levelled) shaped 3-2 Pachner move. Two (levelled) shaped pseudo 3-manifolds $X$ and $Y$ are called \emph{equivalent} if there exist gauge equivalent (levelled) shaped pseudo 3-manifolds $X'$ and $Y'$ which are respective  Pachner refinements of $X$ and $Y$.
\end{definition}

In its original formulation, our Teichm\"{u}ler TQFT is not defined on all levelled shaped pseudo 3-manifolds. It is only guaranteed to be well defined on $S_r(X)$ (since we need the positivity of the angles to make certain integrals absolutely convergent) and when $H_2(X-\Delta_0(X), \bZ) = 0$. The latter condition guarantees that we can multiply the distributions for all the tetrahedra and peform the necessary push forward of this product. We therefore need the following definition.

\begin{definition}\label{admis}
An oriented triangulated pseudo $3$-manifold is called \emph{admissible} if
$S_r(X) \neq \emptyset$
and
$
H_2(X-\Delta_0(X), \bZ) = 0.
$
\end{definition}

The equivalence of admissible levelled shaped pseudo 3-manifolds also needs to be such that all involved pseudo 3-manifolds are admissible, hence we introduce a stronger notion of {\em admissibly equivalence}.

\begin{definition}\label{admiseq}
 Two admissible (levelled) shaped pseudo 3-manifolds $X$ and $Y$ are called \emph{admissibly equivalent} if there exists a gauge equivalence $h\colon X' \ra Y'$ of (levelled) shaped pseudo 3-manifolds $X'$ and $Y'$ which are respective Pachner refinements of $X$ and $Y$, such that 
 $$ \Delta_1(X') = \Delta_1(X) \cup D_X, \mbox{  } \Delta_1(Y') = \Delta_1(Y) \cup D_Y$$
 and
 $$ h(S_r(X') \cap \tilde\Omega_{X',r}(D_X)^{-1}(2\pi)) \cap  \tilde\Omega_{Y',r}(D_Y)^{-1}(2\pi) \neq \emptyset.$$
 The corresponding equivalence classes are called \emph{admissible equivalence classes}.
\end{definition}

\begin{theorem}[\cite{AK1}]\label{equivalencetheo}
Suppose two (levelled) shaped pseudo $3$-manifolds $X$ and $Y$ are equivalent. Then there exist $D\subset \Delta_1(X)$ and $D' \subset \Delta_1(Y)$, a bijection
$$i \colon \Delta_1(X)\setminus D \ra  \Delta_1(Y)\setminus D',$$
and a Poisson isomorphism
$$R\colon  \tilde \Omega_{X, r} (D)^{-1}(2\pi) \ra \tilde \Omega_{Y, r} (D')^{-1}(2\pi),$$
covered by  an affine $\bR$-bundle isomorphism 
$$
\tilde R\colon \widetilde{LS}_r(X)|_{\tilde \Omega_{X, r} (D)^{-1}(2\pi) }\to \widetilde{LS}_r(Y)|_{\tilde \Omega_{Y, r} (D')^{-1}(2\pi)},
$$
such that the following  diagram is commutative
 $$\begin{CD}
\tilde \Omega_{X, r} (D)^{-1}(2\pi)  @>R>>  \tilde \Omega_{Y, r} (D')^{-1}(2\pi)\\
@VV{\operatorname{proj}\ \circ \ \tilde \Omega_{X,r}}V @VV{\operatorname{proj}\ \circ \ \tilde \Omega_{Y,r}}V\\
\bR^{\Delta_1(X)\setminus D} @> i^* >>  \bR^{\Delta_1(Y)\setminus D'}.
\end{CD}
$$
Moreover, if $X$ and $Y$ are admissible and admissibly equivalent, then the isomorphism $R$ takes an open non-empty convex subset $U$ of $S_r(X)\cap \tilde \Omega_{X, r} (D)^{-1}(2\pi) $ onto an open non-empty convex subset $U'$ of  $S_r(Y)\cap \tilde \Omega_{Y, r} (D)^{-1}(2\pi) $. 
\end{theorem}

We observe that in the notation of Definition \ref{admiseq}
$$ D = \Delta_1(X) \cap h^{-1}(D_Y), \mbox{ } D' = \Delta_1(Y) \cap h(D_X).$$

Let us now recall the categroid of admissible levelled shaped pseudo $3$-manifolds.
To this end we first need to recall the underlying category $\CC$ where equivalence classes of levelled shaped pseudo 3-manifolds form morphisms, the objects are triangulated surfaces, and composition is given by gluings along the relevant parts of the boundaries by edge orientation preserving and face orientation reversing CW-homeomorphisms with the obvious composition of dihedral angles and addition of levels. Depending on the way we split the boundary into a source and a target, one and the same levelled shaped pseudo 3-manifold can be interpreted as different morphisms in $\CC$. Nonetheless, there is one canonical choice defined as follows.

For a tetrahedron $T=[v_0,v_1,v_2,v_3]$ in $\REALS^3$ with ordered vertices $v_0,v_1,v_2,v_3$, we define its sign by
\[
\sign(T)=\sign(\det(v_1-v_0,v_2-v_0,v_3-v_0)),
\]
as well as the signs of its faces
\[
\sign(\partial_iT)=(-1)^{i}\sign(T),\quad i\in\{0,\ldots,3\}.
\]
\begin{figure}[htbp]
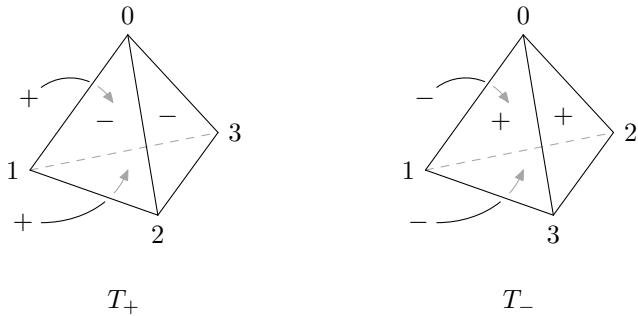

\centering
\includegraphics{tetra-1.mps}
\hskip2cm
\includegraphics{tetra-2.mps}
\caption{Face orientations}
\label{fig:+/-}
\end{figure}
For a pseudo 3-manifold $X$, the signs of the faces of  the tetrahedra of $X$ induce a sign function on the faces of the boundary of $X$,
\[
\sign_X\colon \Delta_2(\partial X)\to\{\pm1\},
\]
which permits to split the boundary of $X$ into two subsets,
\[
\partial X=\partial_+X\cup\partial_-X,\quad
\Delta_2(\partial_\pm X)=\sign_X^{-1}(\pm1),
\]
composed of equal numbers of triangles. For example, in the case of a tetrahedron $T$ with $\sign(T)=1$, we have $\Delta_2(\partial_+T)=\{\partial_0T,\partial_2T\}$, and $\Delta_2(\partial_-T)=\{\partial_1T,\partial_3T\}$.
In what follows, unless specified otherwise, (the equivalence class of) a levelled shaped pseudo 3-manifold $X$  will always be thought of as a $\CC$-morphism between the objects $\partial_-X$ and $\partial_+X$, i.e.
\[
X\in\Hom_\CC(\partial_-X,\partial_+X).
\]
We will also consider more general morphisms to be included in $\CC$, namely morphisms as above, but where we add an allowed one dimensional sub-complexes as defined above. This means that we allow objects where we have special marked vertices on the boundaries and when we compose such morphisms, we assume that all univalent vertices of the sub-complexes, which are contained in the surfaces we glue on, match up pairwise, thus resulting in an allowed one dimensional sub-complex in the morphism obtained by gluing. In the equivalence relation, this one dimensional sub-complex should be carried all the way through the equivalences specified in Definition \ref{PMove} and \ref{admiseq}, but in such a way that we never perform any 3-2 Pachner moves on edges, which are part of the one-dimensional sub-complex. In the rest of this paper we use the term \emph{levelled shaped pseudo $3$-manifold} to mean any morphism of $\CC$ (including also the morphisms we just added to $\CC$).

Our TQFT is not defined on the full category $\CC$, but only on the sub-categroid of admissible equivalence classes of admissible morphisms.

\begin{definition}
The categroid $\CC_a$ of admissible levelled shaped pseudo $3$-manifolds is the sub-categroid of the category of levelled shaped pseudo $3$-manifolds whose morphisms consist of admissible equivalence classes of admissible levelled shaped pseudo $3$-manifolds.
\end{definition}

Gluing in this sub-categroid is the one induced from the category $\CC$ and it is only defined for those pairs of admissible morphisms for which the glued morphism in $\CC$ is also admissible.

There is a concise graphical presentation of pseudo $3$-manifolds introduced in \cite{AK1}. To each tetrahedron $T$, it associates the graph 
 \begin{equation}
   T=\begin{tikzpicture}[yscale=0.4,baseline=5]
     \draw[very thick] (0,0)--(3,0);
      \draw[gray] (0,0)--(0,1) (1,0)--(1,1)(2,0)--(2,1)(3,0)--(3,1);
     \draw (0,1) node[above] {$\partial_0T$};\draw (1,1) node[above]{$\partial_1T$};\draw (2,1) node[above]{$\partial_2T$};\draw (3,1) node[above]{ $\partial_3T$};
   \end{tikzpicture}
 \end{equation}
 where each of the four codimension one faces of $T$ corresponds to a vertical half edge. We connect the half edges according to the face identifications of the tetrahedra in a given triangulated pseudo $3$-manifold (in the case of non-empty boundary, the resulting graph will  also have open half edges).  For example, the graphs for the pseudo $3$-manifolds representing the complements of the trefoil knot $3_1$, the figure eight knot $4_1$ and the $5_2$-knot are as follows
 
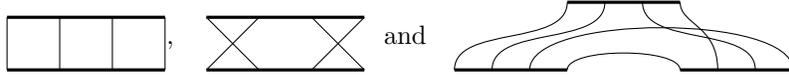
\begin{equation*}\label{P:trefoil}
\begin{tikzpicture}[baseline=10,scale=.7]
\draw[very thick] (0,0)--(3,0);\draw[very thick] (0,1)--(3,1);
\draw (0,0)--(0,1);
\draw (1,0)--(1,1);
\draw (2,0)--(2,1);
\draw (3,0)--(3,1);
\end{tikzpicture},
\ \ \ 
\begin{tikzpicture}[baseline=10,scale=.7]
\draw[very thick] (0,0)--(3,0);\draw[very thick] (0,1)--(3,1);
\draw(0,0)--(1,1);
\draw (1,0)--(0,1);
\draw (2,0)--(3,1);
\draw (3,0)--(2,1);
\end{tikzpicture}\ \ \text{and}
\ \ \ 
 \begin{tikzpicture}[baseline=10,xscale=.5,yscale=.3]
 \draw[very thick] (0,0)--(3,0);
 \draw[very thick] (6,0)--(9,0);
 \draw[very thick] (3,3)--(6,3);
 \draw(3,0)..controls (3,1) and (6,1)..(6,0);
 \draw(0,0)..controls (0,2) and (3,1)..(3,3);
 \draw(1,0)..controls (1,2) and (4,1)..(4,3);
 \draw(2,0)..controls (2,2.5) and (9,2.5)..(9,0);
 \draw(5,3)..controls (5,1) and (8,2)..(8,0);
 \draw(6,3)..controls (6,2) and (7,2)..(7,0);
 \end{tikzpicture}
\end{equation*}
All these examples correspond to 3-manifolds with one cusp, i.e. they are 1-vertex triangulations with the vertex having a neighborhood  homeomorphic to the cone over the torus. As it will be seen below, our TQFT functor, up to overall orientation, can be written down just based on such graphical presentation.

\section{The target categroid}

\label{target}


The target categroid for the Teichm\"{u}ller TQFT is given by  tempered distributions. They form only a categroid since the kind of composition of distributions we have in mind is not defined for all tempered distributions.

Recall that the space of (complex) tempered distributions $\S'(\bR^n)$ is the space of continuous linear functionals on the (complex) Schwartz space $\S(\bR^n)$. By the Schwartz presentation theorem (see e.g. Theorem V.10 p. 139 \cite{RS1}), any tempered distribution can be represented by a finite derivative of a continuous function with polynomial growth, hence we may informally think of tempered distributions as functions defined on $\bR^n$. 
The integral formula 
$$
\varphi(f) = \int_{\bR^n} \varphi(x) f(x) dx.
$$
exhibits the inclusion $\S(\bR^n)\subset \S'(\bR^n)$.

\begin{definition}\label{tcat}
The categroid $\D$ has as objects finite sets and for two finite sets $n,m$ the set of morphisms from $n$ to $m$ is
$$
\Hom_{\D}(n,m) = \S'(\bR^{n\sqcup m}).
$$
\end{definition} 
Denoting by  $ \L( \S(\bR^{n}), \S'(\bR^{m}))$ the space of continuous linear maps from $\S(\bR^{n})$ to $ \S'(\bR^{m})$, we remark that we have an isomorphism
\begin{equation*}
\tilde{\cdot} :   \L( \S(\bR^{n}), \S'(\bR^{m})) \ra \S'(\bR^{n\sqcup m})
\end{equation*}
determined by the formula
$$
\varphi(f)(g) = \tilde{\varphi}(f\otimes g)
$$
for all $\varphi \in  \L( \S(\bR^{n}), \S'(\bR^{m}))$, $f\in  \S(\bR^{n})$, and $g\in \S(\bR^{m})$. This is the content of the Nuclear theorem, see e.g. \cite{RS1}, Theorem V.12, p. 141. The reason why we get a categroid rather than a category is because we cannot compose all composable (in the usual categorical sense) morphisms, but only a subset thereof.
The partially defined composition in this categroid is  defined as follows. 

Let $n,m,l$ be three finite sets, $A\in \Hom_{\D}(n,m)$ and $B\in \Hom_{\D}(m,l)$. According to the tempered distribution analog of Theorem 6.1.2. in \cite{Hor1}, associated to the projections 
$$
\pi_{n,m}\colon \bR^{n\sqcup m \sqcup l}\to \bR^{n\sqcup m},\quad \pi_{m,l}\colon \bR^{n\sqcup m \sqcup l}\to \bR^{m\sqcup l},
$$
we have the pull back maps
$$
\pi_{n,m}^* \colon \S'(\bR^{n\sqcup m}) \ra \S'(\bR^{n\sqcup m \sqcup l}) \mbox{ and } \pi_{m,l}^* \colon\S'(\bR^{m\sqcup l}) \ra \S'(\bR^{n\sqcup m \sqcup l}).
$$
By theorem IX.45 in \cite{RS2} (see also Appendix~B in \cite{AK1}), the product
$$
\pi_{n,m}^*(A)\pi_{m,l}^*(B) \in \S'(\bR^{n\sqcup m \sqcup l})
$$
is well defined provided the wave front sets of $\pi_{n,m}^*(A)$ and $\pi_{m,l}^*(B)$ satisfy the following transversality condition
\begin{equation}\label{wftrans}
(\wf(\pi_{n,m}^*(A)) \oplus \wf(\pi_{m,l}^*(B)) )\cap Z_{n\sqcup m \sqcup l} = \varnothing,
\end{equation}
where $Z_{n\sqcup m \sqcup l}$ is the zero section of $T^*(\bR^{n\sqcup m \sqcup l})$. If we now further assume that $\pi_{n,m}^*(A)\pi_{m,l}^*(B)$ continuously extends to $\S(\bR^{n\sqcup m \sqcup l})_m$ as is defined in Appendix~B of \cite{AK1},
then we obtain a well defined element
$$
(\pi_{n,l})_*(\pi_{n,m}^*(A)\pi_{m,l}^*(B) ) \in \S'(\bR^{n\sqcup l}).
$$
\begin{definition}
For $A\in \Hom_{\D}(n,m)$ and $B\in \Hom_{\D}(m,l)$ satisfying condition~(\ref{wftrans}) and such that $\pi_{n,m}^*(A)\pi_{m,l}^*(B)$ continuously extends to a well defined element of the dual of $S(\bR^{n\sqcup m \sqcup l})_m$, we define
$$
AB = (\pi_{n,l})_*(\pi_{n,m}^*(A)\pi_{m,l}^*(B) ) \in \Hom_{\D}(n,l).
$$
\end{definition}
For any $A\in  \L( \S(\bR^{n}), \S'(\bR^{m}))$, we have unique adjoint $A^*\in \L( \S(\bR^{m}), \S'(\bR^{n}))$ defined by the formula
$$
A^*(f)(g) = \overline{A(\bar{g})(\bar f)}
$$ for all $f\in  \S(\bR^{m})$ and $g\in \S(\bR^{n})$.

\section{The TQFT functor}

\label{functor}

We shall describe the Teichm\"{u}ller TQFT functor $F_\hbar$ from \cite{AK1}. 
First we recall the definition of a $*$-functor in our context.

\begin{definition}
A functor $F\colon \CC_a \to \D$ is said to be a  \emph{$*$-functor} if
\[
F(X^*)=F(X)^*,
\]
where $X^*$ is $X$ with opposite orientation, and $F(X)^*$ is the adjoint of $F(X)$.
\end{definition}

On the level of objects we define
$$ F_\hbar(\Sigma) = \Delta_2(\Sigma),\quad \forall\Sigma\in\operatorname{Ob}\mathcal{B}_a.$$

In order to define $ F_\hbar$ on morphisms, we need a special function called Faddeev's quantum dilogarithm \cite{F}.

\begin{definition} \emph{Faddeev's quantum dilogarithm} is the function of two complex arguments $z$ and $\la$ defined for $|\im z| <\frac12 |\la+\la^{-1}|$  by the formula 
\[
\QDILOG(z):=\exp\left(
\int_{C}
\frac{e^{-2\IMUN zw}\, dw}{4\sinh(w\la)
\sinh(w/\la) w}\right),
\]
where the contour $C$ runs along the real axis deviating into the upper half plane in the vicinity of the origin, and extended by the functional equation
\begin{equation*}
\QDILOG(z-\IMUN\la^{\pm1}/2)=(1+e^{2\pi\la^{\pm1}z })
\QDILOG(z+\IMUN\la^{\pm1}/2)
\end{equation*}
to a meromorphic function for $z\in \bC$.
\end{definition}

It is easily seen that  $\QDILOG(z)$ depends on $\la$ only through the combination $\hbar$ defined by the formula
\[
\hbar:=\left(\la+\la^{-1}\right)^{-2}.
\]
In what follows, we assume that the complex parameter $\la$ is such that $\hbar\in\bR_{>0}$. This assumption corresponds to a \emph{unitary}  TQFT, but, in case of need, one can easily go to arbitrary $\la\in \bC\setminus \imun\bR$ by analytic continuation.

The value of $F_\hbar$ on a morphism  $(X,\Gamma)$ of $\CC_a$ is given by the formula which singles out the dependence on the level
$$F_\hbar(X,\Gamma) = e^{\imun\pi\frac{\ell_X}{4\hbar}}\parfun(X,\Gamma)\in\S'\left(\bR^{\Delta_2(\partial X)}\right),$$
where $Z_\hbar(X,\Gamma) $
is the level independent part.

The value of $Z_\hbar$ on the morphism  of $\CC_a$ determined by a single tetrahedron $T$ with $\sign(T)=1$ is an element
$$ Z_\hbar(T,\alpha_T) \in \S'(\bR^{\Delta_2(T)})$$
given by the explicit formula
\begin{equation}\label{tet-int-f}
Z_\hbar(T,\alpha_T)(x_0,x_1,x_2, x_3)=\delta(x_0-x_1+x_2)\frac{e^{2\pi\imun (x_3-x_2)\left(x_0+\frac{\alpha_3}{2\imun\sqrt{\hbar}}\right)+\pi \imun\frac{\varphi_T}{4\hbar}}}{\QDILOG\left(x_3-x_2+\frac{1-\alpha_1}{2\imun\sqrt{\hbar}}\right)}
\end{equation}
where $\delta$ is Dirac's delta-function supported at $0\in\bR$,
\[
\varphi_T:=\alpha_1\alpha_3+\frac{\alpha_1-\alpha_3}3-\frac{2\hbar+1}{6},\quad
\alpha_i:=\frac1\pi\alpha_T(\partial_0\partial_iT),\quad i\in\{1,2,3\},
\]
and
\[
x_i:=x(\partial_iT),\quad x\colon \Delta_2(\partial T)\to\REALS.
\]
For a negative tetrahedron $\bar T$ with $\sign(\bar T)=-1$ we set
$$Z_\hbar(\bar T) = Z_\hbar(T)^*.$$
It is not hard to check that these assignments give tempered distributions provided $\alpha_i>0$, $i=1,2,3$. 

The value of $Z_\hbar$ on arbitrary morphism $(X,\Gamma)$ in $\CC_a$ is given by composing all the distributions $Z_\hbar(T)$, where $T$ runs over $\Delta_3(X)$, according to the face identifications which build $X$ out of the disjoint union
$$ \tilde X:= \bigsqcup_{T\in \Delta_3(X)} T.$$
By using the graphical presentation of $X$ described above, with the additional information on the orientation of $X$, the prescription is as follows. One should label the thin edges of the graph with variables $x_i$, where $i=1, \ldots,|\Delta_2(X)|$, then take the product over all tetrahedra of the expression (\ref{tet-int-f}) or its complex conjugate adapted to each tetrahedron in accordance with the variables attached to its four faces and the dihedral angles, and integrate over all real values of  $x_i$.

Let us illustrate this with the example of the complement $X$ of the knot $5_2$ which is represented by the diagram
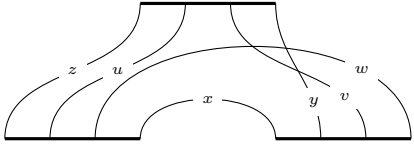
\[
 \begin{tikzpicture}[scale=.6]
 \draw[very thick] (0,0)--(3,0);
 \draw[very thick] (6,0)--(9,0);
 \draw[very thick] (3,3)--(6,3);
 \draw (3,0) to [out=90, in=90, edge node={node[fill=white] {\tiny $x$}}]
 (6,0);
 \draw(0,0)to [out=90, in=-90, edge node={node[fill=white] {\tiny $z$}}]
 (3,3);
 \draw(1,0)to [out=90, in=-90, edge node={node[fill=white] {\tiny $u$}}]
 (4,3);
 \draw(2,0)to [out=90, in=90, edge node={node[fill=white, near end] {\tiny $w$}}]
 (9,0);
 \draw(5,3)to [out=-90, in=90, edge node={node[fill=white, near end] {\tiny $v$}}]
 (8,0);
 \draw(6,3)to [out=-90, in=90, edge node={node[fill=white, near end] {\tiny $y$}}]
 (7,0);
 \end{tikzpicture}
\]
We denote $T_1,T_2,T_3$ the left, right, and top tetrahedra respectively with their dihedral angles $\alpha_{T_i} = 2\pi(a_i,b_i,c_i)$, such that $a_i+b_i+c_i = \frac12$, $i=1,2,3$. We choose the orientation so that  all tetrahedra are positive, and
we impose the conditions that all edges are balanced ($\Gamma = \emptyset$) which correspond to two equations
\[
2a_3=a_1+c_2,\quad b_3=c_1+b_2.
\]
We thus get by the definition of our TQFT that
\begin{multline}\label{52c}
Z_\hbar(X)=\int_{\REALS^6}
Z_\hbar(T_1,\alpha_{T_1})(z,u,w,x) Z_\hbar(T_2,\alpha_{T_2})(x,y,v,w) \\
\times Z_\hbar(T_3,\alpha_{T_3})(y,v,u,z)\operatorname{d}^6(x,y,z,u,v,w),
\end{multline}
where we observe that the integrand indeed extends to the dual of $S(\REALS^6)_6$, thus it can be pushed forward to a point, which is the precise meaning of the integral in (\ref{52c}).
The calculation in Section 11.6 of \cite{AK1} gives 
$$Z_\hbar(X) = \nu_{c_1,b_1}\nu_{b_2,a_2}\nu_{c_3,b_3}e^{\imun\pi\cla^2(1-2a_1)(1-2c_2)}
\int_{2\cla (a_1-a_3)+\REALS}\chi_{5_2}(x,\lambda) \operatorname{d}\!x,$$
where $ \lambda := a_1-c_1+b_2-a_3$,
\[
\nu_{a,b}: = e^{4\pi \imun \cla^2a(a+b)} e^{-\pi \imun \cla^2(4(a-b) +1)/6},\quad \cla:=\frac{\imun}{2\sqrt{\hbar}}=\frac{\imun}{2}(\la+\la^{-1}),
\]
and
\begin{multline}\label{eq:chi5_2}
\chi_{5_2}(x,\lambda):=\chi_{5_2}(x)e^{4 \pi \imun\cla x\lambda},\\ 
\chi_{5_2}(x):=e^{-\imun\pi/3}\int_{\REALS-\imun0}\operatorname{d}\!z\, \frac{
e^{\imun\pi(z-x) (z+x)}}{\QDILOG(z+x)\QDILOG(z-x)\QDILOG(z)}.
\end{multline}

\vskip.5cm

Returning now back to the case of a general $X$, we need to know that all the compositions of the $Z_\hbar(T,\alpha_T)$'s are allowed in $\D$. This is precisely the content of Theorem 9 in \cite{AK1}, which establishes that for admissible $X$, the wave front sets of the distributions $\pi_T^* Z_\hbar(T)$,  where $\pi_T : \bR^{\Delta_2(X)} \ra \bR^{\Delta_2(T)}$ is the natural projection for each $T\in \Delta_3(X)$, are transverse and hence they can be multiplied and their product can be pulled back to $\bR^{\Delta_2(X)}$ and pushed forward along the projection from $\bR^{\Delta_2(X)}$ to $\bR^{\Delta_2(\partial X)}$. 

We emphasize that for an admissible pseudo $3$-manifold $X$ together with an allowed sub-complex $\Gamma$ of $\Delta_1(X)$, our TQFT functor provides us with the following well defined \emph{real analytic} function
$$
F_\hbar(X,\Gamma) : LS_r(X)\cap\tilde{\Omega}_{X,r}(E_\Gamma)^{-1}(2\pi) \ra \S'(\bR^{\partial X}),
$$
where $E_\Gamma$ is the set of internal edges of $X$ which are not in $\Gamma$. We note that if $(X,\Gamma)$ is admissibly equivalent to $(X',\Gamma')$, then Theorem \ref{equivalencetheo} provides an explicit affine map from a non-empty open convex subset of $LS_r(X)\cap\tilde{\Omega}_{X,r}(E_\Gamma)^{-1}(2\pi)$ to an open convex subset of $LS_r(X')\cap\tilde{\Omega}_{X',r}(E_{\Gamma'})^{-1}(2\pi)$ and under this map the restrictions of $F_\hbar(X,\Gamma)$ and $F_\hbar(X',\Gamma')$ to these two non-empty convex open subsets agree. Since a real analytic map defined on an open convex set of some Euclidian space is uniquely determined by its restriction to any smaller non-empty open convex subset, we see that $F_\hbar(X,\Gamma)$  and  $F_\hbar(X',\Gamma')$  uniquely determine each other on their domains of definition. It is in this sense that our Teichm\"{u}ller TQFT is well-defined on the set of equivalence classes of admissible levelled shaped pseudo $3$-manifolds.

For the case $\partial X = \varnothing$, we have $\S'(\bR^{\partial X})= \bC$ and so, in this case, we simply get a complex valued function on $LS_r(X) \cap\tilde{\Omega}_{X,r}(E_\Gamma)^{-1}(2\pi)$.
In particular, the value of the functor $F_\hbar$ on any fully balanced admissible levelled shaped 3-manifold is a complex number, which is a topological invariant, in the sense that if two fully balanced admissible levelled shaped 3-manifolds are admissibly equivalent, then $F_\hbar$ assigns one and the same complex number to them. We recall that fully balanced means that all edges are balanced.

Our main Theorem 4 of \cite{AK1} now guarantees that this assignment, in fact, gives a well-defined functor.

\begin{theorem}\label{Main}\label{MT}
For any $\hbar\in\REALS_{>0}$, the above assignment defines a  $*$-functor
$$F_\hbar\colon \CC_a \to \D$$
which we call the Teichm\"{u}ller TQFT.
\end{theorem}

\section{The volume conjecture for the Teichm\"{u}ller TQFT}
\label{vol}

In this subsection we recall our conjecture from \cite{AK1} concerning our Teichm\"{u}ller TQFT $F_\hbar$, which, among other things, provides a relation to the hyperbolic volume in the asymptotic limit $\hbar\to0$. 

\begin{conjecture}\label{conj}
Let $M$ be a closed oriented compact 3-manifold.
For any hyperbolic knot $K\subset  M$, there exists a smooth function
$J_{M,K}(\hbar,x)$ on $\REALS_{>0}\times\REALS$  which has the following properties.
\begin{enumerate}
\item
For any fully balanced shaped ideal triangulation $X$ of the complement  of $K$ in $M$, there exist a gauge invariant real  linear combination of dihedral angles $\lambda$ and  a (gauge non-invariant) real quadratic polynomial of dihedral angles $\phi$ such that
\[
\parfun(X)=e^{\imun\frac{\phi}{\hbar}}\int_{\REALS} J_{M,K}(\hbar,x)e^{-\frac{x\lambda}{\sqrt{\hbar}}} \operatorname{d}\!x
\]
\item The hyperbolic volume of the complement of $K$ in $M$ is recovered as the following limit
\[
\lim_{\hbar\to0}2\pi\hbar\log\vert J_{M,K}(\hbar,0)\vert=-\vol(M\setminus K).
\]
\end{enumerate}
\end{conjecture}
\begin{remark}
It is very important to notice that we in part (2) of this conjecture have a negative sign in the right hand side, which differs from the volume conjecture of \cite{K6}. In this case, the invariant exponentially decays (rather than grows) with the rate being given by the hyperbolic volume.
\end{remark}

In \cite{AK1}, we checked this conjecture for the first two hyperbolic knots.

\begin{theorem}\label{th:4-1--5-2}
Conjecture \ref{conj} is true for the pairs $(S^3,4_1)$ and $(S^3,5_2)$ with
\[
J_{S^3,4_1}(\hbar,x)=\chi_{4_1}(x),\quad J_{S^3,5_2}(\hbar,x)=\chi_{5_2}(x),
\]
where the functions $\chi_{4_1}(x)$ and $\chi_{5_2}(x)$ are given by
$$ \chi_{4_1}(x) = \int_{\bR-i0} \frac{\Phi_\la(x-y)}{\Phi_\la(y)} e^{2\pi i x(2y-x)} \operatorname{d}\!y$$
and  $\chi_{5_2}(x)$ is given in (\ref{eq:chi5_2}) above.
\end{theorem}

See also \cite{AM} for a precise statement of the generalisation of the above conjecture to the higher level generalisation of the Teichm\"{u}ller TQFT and more examples in \cite{AN}. Further, in \cite{AMa} we have presented the precise formulation of the AJ-conjecture for the Teichm\"{u}ller TQFT.

\section{Future perspectives}
\label{future}

In the paper \cite{AK4} we have presented a new formulation of the Teichm\"{u}ller TQFT and a further higher level generalization of the theory, which we think of as a version of the complex quantum Chern-Simons theory \cite{AK5} (see also \cite{AM}). Let us here briefly recall this new formulation and state some predictions for the further perspectives for the Teichm\"{u}ller TQFT.

The main player behind the new formulation of the Teichm\"{u}ller TQFT is the edge-face tranform using the Weil--Gel'fand--Zak (WGZ) transformation, which we now recall.

We consider the following multiplier construction for a line bundle over the two torus. We have the natural translation action of $\bZ^2$ on $\bR^2$ with the quotient $\Pi := {\mathbb S}^2$ where ${\mathbb S} := \bR/\bZ$.
Consider the following multipliers 
\begin{equation*}
\varphi\colon \bZ^2\times\bR^2\to U(1),\quad\varphi((m,n),(x,y))=(-1)^{mn}e^{\pi\Imun(nx-my)},
\end{equation*}
which induce an action of ${\mathbb Z}^2$ on the trivial bundle $\bR^2\times\bC$ and we define 
$$L = (\bR^2\times\bC)/\bZ^2$$  
as a complex line bundle over $\Pi$. 

We define the Weil--Gel'fand--Zak (WGZ) transformation 
$$ W : \mathcal{S}(\REALS) \ra C^\infty(\Pi, L)$$
by the formula
\begin{equation*}\label{eq:wgz}
(Wf)(x,y)=e^{\pi\Imun xy}\sum_{m\in\INTEGERS}f(x+m)e^{2\pi\Imun my}.
\end{equation*}

\begin{proposition}
The WGZ-transformation 
$$ W : \mathcal{S}(\REALS) \ra C^\infty(\Pi, L)$$
 is an isomorphism of Fr\'{e}chet spaces with the inverse explicitly given by the formula
\[
(W^{-1}g)(x)=\int_0^1g(x,y)e^{-\pi\Imun xy}\operatorname{d}\!y
\]
for any $g\in C^\infty(\Pi, L)$.
\end{proposition}

The same statement also holds for
$$ \overline{W} : \mathcal{S}(\REALS) \ra C^\infty(\Pi, L^*)$$
defined by
$$\overline{W}(f)(x,y)=W (f)(x,-y).$$

Guided by equation (\ref{tet-int-f}), we define for any $f\in {\mathcal S}(\REALS)$ the associated tempered distribution $H(f)\in {\mathcal S}'(\REALS^4)$ as follows
\[
H(f)(x_0,x_1,x_2,x_3) = \delta(x_0+x_2-x_1)f(x_3-x_2)e^{2\pi\Imun x_0(x_3-x_2)},
\]
which we consider as a continuous linear map
$$H(f) : {\mathcal S}(\bR^2) \ra {\mathcal S}'(\bR^2)$$
via the formula
$$H(f)(g)(h) = H(f)(\hat\pi^*(g)\check\pi^*(h)).$$
Here $\hat \pi,\check\pi \colon \bR^4 \ra \bR^2$ are given by $\hat\pi(x_0,x_1,x_2,x_3) = (x_1,x_3)$ and $\check\pi(x_0,x_1,x_2,x_3) = (x_0,x_2)$.
We now consider the following tensor extension of the WGZ-transform 
$$
W\otimes  W : {\mathcal S}(\REALS^2) \ra C^\infty(\Pi \times \Pi, L\boxtimes L)
$$
defined by
\begin{eqnarray*} \lefteqn{W\otimes W (h)(s,t,x,y) = }\\
& \phantom{kkkk}e^{\pi i (st + xy)} \sum_{m_1,m_2\in \bZ} h(s+m_1,x+m_2) e^{2\pi i(m_1t + m_2 y)}.
\end{eqnarray*}
and similarly
$$\overline{W}\otimes\overline{W} : {\mathcal S}(\bR^2) \ra C^\infty(\Pi \times \Pi, L^*\boxtimes L^*),$$
given by
\begin{eqnarray*} \lefteqn{\overline{W}\otimes\overline{W}(h)(s,t,;x,y) = }\\
& \phantom{kkkk}e^{-\pi i (st + xy)} \sum_{m_1,m_2\in \bZ} h(s+m_1,x+m_2) e^{-2\pi i(m_1t + m_2 y)}.
\end{eqnarray*}
Consider now the maps
$$ F : \bS^5 \ra \Pi, \hskip4mm \tilde\pi_i :  \bS^5 \ra \Pi^2,\quad i=1,2,$$
given by
$$ F(u,s,t,x,y) = (u, s+t+u-y),$$
$$\tilde\pi_1(u,s,t,x,y) = (s+x,t+u,x+u,y-t-u)$$
and $\tilde\pi_2$ is the map which projects away the first factor and onto the last four factors. We then get that there exists a natural isomorphism
$$ F^*L \otimes \tilde\pi_1^*(L\boxtimes L) \cong  \tilde\pi_2^*(L\boxtimes L).$$
There is an obvious embedding
$$ C^\infty(\Pi^2, L\boxtimes L) \subset  C^\infty(\Pi^2, L^*\boxtimes L^*)^*,$$
obtained by pointwise evaluation followed by integration over $\Pi^2$. Thus we see that
$$ ( \tilde\pi_2)_* (F^*W(f) \tilde\pi_1^*(W\otimes W)(g)) \in C^\infty(\Pi^2, L^*\boxtimes L^*)^*$$
for any $g\in {\mathcal S}(\bR^2)$.

Let $T : \bS^3 \ra \Pi$ be given by $T(x,y,z) = (y-x,z-y)$ and further $E : \bS^6 \ra \Pi$ by
$$ E(x_{01},x_{02},x_{03},x_{12},x_{13},x_{23})= (x_{02}+x_{13}-x_{03}-x_{21},x_{02}+x_{13}-x_{01}-x_{23}).$$
We also introduce the following two maps $\pi_i : \bS^6 \ra \Pi^2$ , $i=1,2$, given by
$$\pi_1(x_{01},x_{02},x_{03},x_{12},x_{13},x_{23}) = T\times T (x_{23},x_{03},x_{02},x_{12},x_{02}.x_{01})$$
and
$$\pi_2(x_{01},x_{02},x_{03},x_{12},x_{13},x_{23}) = T\times T (x_{23},x_{13},x_{12},x_{13},x_{03},x_{01}).$$
Finally we further need the map $P : \bS^6 \ra \bS^5$ given by
\begin{eqnarray*}
& & P(x_{01},x_{02},x_{03},x_{12},x_{13},x_{23}) \\
& &\phantom{mmmmm} = ( x_{02}+x_{13}-x_{03}-x_{12}, x_{13} - x_{23}, x_{12} - x_{13}, x_{03} - x_{13}, x_{01} - x_{03})
\end{eqnarray*}
It is elementary to verify that 
$$ E = F \circ P  \text { and } \pi_i = \tilde\pi_i \circ P.$$
We now arrive at the important ``edge-face transformation" as established in \cite{AK4}.
\begin{theorem}\label{eft}
The distribution $$H(f)\in {\mathcal S}'(\bR^4)$$ and the section $$E^*W(f)\in C^\infty(\bS^6, E^*L)$$ are related by the formula
$$
H(f)(\hat\pi^*(g)\check\pi^*(h)) = \pi^{(6)}_* (E^*W(f)\pi_1^*(W\otimes W)(g) \pi_2^*(\overline{W}\otimes\overline{W})(h))
$$
for all $g,h\in S(\bR^2)$.
\end{theorem}

This theorem is the corner stone in understanding how to transform our original Teichm\"{u}ller TQFT, which is build from the fundamental distribution 
$$F_\hbar(T) \in {\mathcal S}'(\bR^{\Delta_2(T)}),$$ 
to our new formulation ${\mathcal F}_\hbar$ of the Teichm\"{u}ller TQFT, where the state variables live on the edges instead and the fundamental object associated to a tetrahedron is the section
$$ {\mathcal F}_\hbar(T) \in C^\infty(\bS^{\Delta_1(T)}, E^*L).$$ 
Here we implicitly use the following notation
 \begin{equation*}
x_{ij}:= x(v_T(i)v_T(j))
\end{equation*}
for any $x\in \bS^{\Delta_1(T)}$. Now, Theorem \ref{eft} simply tells us that the two distributions $F_\hbar(T)$ and ${\mathcal F}_\hbar(T)$ are related via the tensor square of the WGZ-transform $W$.

Let us recall the formula for ${\mathcal F}_\hbar(T)$ from \cite{AK4}
$$ {\mathcal F}_\hbar(T) = E^*(g_{\alpha_0,\alpha_2}),$$
where, for two positive real numbers $a$ and $c$ satisfying $a+c<1/2$, we  let
\begin{equation*}\label{eq:gac}
 g_{a,c} = W(\tilde\psi'_{a,c}),
\end{equation*}
\begin{equation*}\label{eq:tpsip}
\tilde\psi'_{a,c}(s):= e^{-\pi \Imun s^2}\tilde\psi_{a,c}(s),
\end{equation*}
\begin{equation*}\label{eq:tpsi}
\tilde\psi_{a,c}(s):= \int_\REALS\psi_{a,c}(t)e^{-2\pi \Imun st}dt,
\end{equation*}
and
\begin{equation*}\label{eq:psi}
\psi_{a,c}(t):= \BQDILOG(t-2\cla(a+c))e^{-4\pi \Imun\cla a(t-\cla(a+c))}
e^{-\pi\Imun\cla^2(4(a-c)+1)/6},
\end{equation*}
where we use the notation $\BQDILOG(x):= 1/\QDILOG(x)$.

As it is described in \cite{AK4} we get a new formulation ${\mathcal F}_\hbar$ of the Teichm\"{u}ller TQFT following the same lines as discussed above determining ${\mathcal F}_\hbar$ on all objects of $\CC_a$ by a similar gluing construction. 

However, this new formulation has several advantages. It allows us to extend the partition function to complex dihedral angles and as such it depends meromorphically on these complexified angles. This allows us to actually establish that we do not need the condition of admissibility and that this new formulation ${\mathcal F}_\hbar$ is well-defined on the full bordism category $\CC$ consisting of equivalence classes of levelled shaped pseudo 3-manifolds. We stress that this means that the functor ${\mathcal F}_\hbar$ is in fact invariant under all 2-3 and 3-2 Pachner moves. 

By modifying the construction of the functor ${\mathcal F}_\hbar$, we can further extend this functor to a version which  depends on a first cohomology class with coefficient in $\bS$ (see \cite{AK4}) and this generalised new version can be related to the original functor $F_\hbar$ by integration over this first cohomology group.

We can use similar ideas (see also \cite{AK4}) to produce a meromorphic extension of $F_\hbar$ to complex angles and as such we can  establish that this theory is also invariant under all 2-3 and 3-2 Pachner moves as argued in \cite{AK4}. Furthermore, via a certain gauge fixing technique also described in \cite{AK4}, we can extend the original theory $F_\hbar$ to be defined on the full bordism category $\CC$, and we describe the precise relation between this extension of the original formulation and the new one in \cite{AK4}.

\end{document}